\pdfoutput=1
\documentclass[a4paper,10pt]{amsart}
\usepackage{amsmath}
\usepackage{amssymb}
\usepackage{amsxtra}
\usepackage{amsthm,amsfonts,latexsym}
\usepackage{epic}
\usepackage{eepic}
\usepackage{fancybox}
\usepackage{boxedminipage}
\usepackage{calc}
\usepackage{tabularx}
\usepackage{graphicx,color}
\usepackage{longtable,booktabs}
\usepackage{supertabular,longtable}
\usepackage{paralist}
\usepackage{amsmath}
\usepackage{amsfonts}
\usepackage{amssymb}
\usepackage{amsthm}
\usepackage[english]{babel}
\usepackage[latin1,ansinew]{inputenc}
\usepackage{a4wide}
\usepackage{pdfpages}

\usepackage{color}

\newcommand{\bea}{\begin{eqnarray}}
\newcommand{\eea}{\end{eqnarray}}
\newcommand{\beann}{\begin{eqnarray*}}
\newcommand{\eeann}{\end{eqnarray*}}
\newcommand{\ba}{\begin{array}}
\newcommand{\ea}{\end{array}}

\newcommand{\beq}{\begin{equation}}
\newcommand{\eeq}{\end{equation}}
\newcommand{\be}{\begin{equation}}
\newcommand{\ee}{\end{equation}}

\accentedsymbol{\hcirc}{ {\overset{\scriptscriptstyle \circ }{ {\rm H} }}}

\newcommand{\R}{\mathbb{R}}

\newcommand{\divv}{ \nabla \!\! \cdot \!}
\newcommand{\bu}{{\bf u}}

\newcommand{\ovrho}{{\overline\rho}}
\newcommand{\unrho}{{\underline\rho}}
\newcommand{\ovc}{{\overline c}}
\newcommand{\unc}{{\underline c}}

\newtheorem{theorem}{Theorem}
\title{Diffuse planar phase boundaries in a two-phase fluid\\ with one incompressible phase}
\author{Heinrich Freist\"uhler and Matthias Kotschote}
\date{June 1, 2013}
\begin{document}
\begin{abstract}
This note studies a family of Navier-Stokes-Allen-Cahn systems
parameterized by temperature. Derived from an internal energy
that corresponds to one incompressible and one compressible phase,   
this family is considered as a simple model for water. 
Decreasing temperature across a critical value, 
a transition takes places from a situation without 
towards one with planar diffuse phase boundaries.
\end{abstract}
\maketitle
\parindent=0cm
In this note, we consider the 
Navier-Stokes-Allen-Cahn system
\begin{alignat}{3}
\partial_t \rho + \divv ( \rho \bu ) & = 0, & \quad 
\notag
\\
\partial_t (\rho \bu)  + \divv ( \rho \bu \otimes\bu + p(\rho,c){\bf I}) 
& = 
\divv \left(
\mu (\nabla\bu + (\nabla \bu)^T)
+
(\lambda \divv \bu) {\bf I} 
- \delta \rho \nabla c \otimes \nabla c 
\right), 
\label{NSAC}
\\
\partial_t (\rho c)  + \divv (\rho c\, \bu ) & =
\delta^{-1/2}(\rho q(\rho,c)
+
\divv \left(\delta \rho \nabla c \right))
\notag
\end{alignat}
of evolutionary partial differential equations that model the spatiotemporal 
behaviour of a compressible viscous or inviscid fluid. 
The fluid is assumed to have a constant temperature $\theta>0$ and to
be a locally homogeneous mixture of two components such that
its local state is completely described by the mass fraction $c$ of one of the
components and the mass, per volume, of the mixture, $\rho$. This density $\rho$
is the reciprocal value,
$$
\rho=1/\tau,
$$
of the fluid's specific volume  $\tau$. 
The behaviour of the fluid is described by a thermodynamic potential 
$$
\bar U(\tau,c,|\nabla c|)\
=
U(\tau,c)
+
{1\over 2}\delta|\nabla c|^2,
\quad
U(\tau,c)
=
\hat U(\tau,c) + W(c,\theta),
$$ 
in which 
$W(c,\theta)$
is the mixing energy and $\delta$ a positive constant.
The pressure $p$ and the transformation rate $q$ derive
from the potential as 
\begin{equation*}
p(\rho,c)= \tilde p(\tau,c) = - U_\tau(\tau,c),
\quad\quad
q(\rho,c)= \tilde q(\tau,c) = - U_c(\tau,c).
\end{equation*}
\par\medskip
System (1) was derived by Blesgen \cite{Ble} and has recently been shown by 
Kotschote  to possess strong solutions \cite{MK3}.  \\

The following two theorems have been proven in \cite{F3} under certain 
assumptions on $\hat U$ and $W$.

\begin{theorem}(Maxwell states and no-flux phase boundaries.)
With $\tilde \theta < \theta_*$ sufficiently close to a critical temperature
$\theta_*$, the following holds for every $\theta\in(\tilde\theta,\theta_*]$. 
There are 
locally uniquely determined fluid states 
$
(\unrho_0,\unc_0),(\ovrho_0,\ovc_0),
$
depending continuously on $\theta$, such that (i)
\begin{align*}
q(\unrho_0,\unc_0) & = 
q(\ovrho_0,\ovc_0)=0, \\
p(\unrho_0,\unc_0) & = 
p(\ovrho_0,\ovc_0) 
\end{align*}
with
$$
\quad\quad\quad
(\unrho_0,\unc_0) = 
(\ovrho_0,\ovc_0)
\ \ \  
\text{if }\ \  \theta=\theta_*, 
$$
and 
(ii) 
if $\theta<\theta_*$, then
$$
\unrho_0 
<\ovrho_0 
$$
and system \eqref{NSAC} admits a no-flux ($m=0$) 
phase boundary 
$$
(\overrightarrow\rho(x),0,\overrightarrow c(x))
\quad\hbox{with}\quad 
(\overrightarrow\rho(-\infty),\overrightarrow 
c(-\infty))=(\unrho_0,\unc_0),
\ \
(\overrightarrow\rho(\infty),\overrightarrow 
c(\infty))=(\ovrho_0,\ovc_0)
$$ 
and (equivalently via $x\mapsto -x$)
a no-flux phase boundary 
$$
(\overleftarrow\rho(x),0,\overleftarrow c(x))
\quad\hbox{with}\quad 
(\overleftarrow\rho(-\infty),\overleftarrow 
c(-\infty))=(\ovrho_0,\ovc_0),
\ \
(\overleftarrow\rho(\infty),\overleftarrow 
c(\infty))=(\unrho_0,\unc_0).
$$ 
\end{theorem}

\begin{theorem}
For sufficiently small mass flux
$
m\neq 0,
$
\par\medskip
(i) 
the (left endstate, profile, right endstate) triple 
$$
(\unrho_0,0,\unc_0),
(\overrightarrow\rho,0,\overrightarrow c),
(\ovrho_0,0,\ovc_0)
$$
perturbs regularly to a (left endstate, profile, right endstate) 
triple 
$$
(\overrightarrow\rho_m^-,\overrightarrow u_m^-,\overrightarrow c_m^-),
(\overrightarrow\rho_{m},\overrightarrow u_m,  \overrightarrow c_{m}),
(\overrightarrow\rho_m^+,\overrightarrow u_m^+,\overrightarrow c_m^+),
$$
corresponding to a traveling-wave phase boundary 
that is densifying if $m>0$ and rarefying if $m<0$; 
\par\medskip
(ii) 
the 
(left endstate, profile, right endstate, profile) 
triple 
$$
(\ovrho_0,0,\ovc_0),
(\overleftarrow\rho,0,\overleftarrow c),
(\unrho_0,0,\unc_0)
$$
perturbs regularly to a (left endstate, profile, right endstate) 
triple 
$$
(\overleftarrow\rho_m^-,\overleftarrow u_m^-,\overleftarrow c_m^-),
(\overleftarrow\rho_m,  \overleftarrow u_m,  \overleftarrow c_m),
(\overleftarrow\rho_m^+,\overleftarrow u_m^+,\overleftarrow c_m^+)
$$
corresponding to a traveling-wave phase boundary 
that is rarifying if $m>0$ and densifying if $m<0$. 
\end{theorem}

The present note serves to point out that Theorems 1 and 2 also hold 
under the following \\

{\bf Modelling Assumptions.} {\it (i) $\hat U$ is of the form
$$
\hat U(\tau,c)=-(1-c)\log\frac{\tau-c\tau_1}{1-c}, 
$$
where $\tau_1$ is a fixed value with
$$
0<\tau_1<1
$$ 
and $\tau$ and $c$ range as 
$$
\tau_1<\tau<\infty
\quad\hbox{and}\quad  
0<c<1.
$$
(ii) With certain critical parameter values $c_*\in(0,1),\theta_*\in\R$, $$W(.,\theta)$$
undergoes a generic transition from convex for $\theta>\theta_*$ to convex-concave-convex
(``double-well'')
for $\theta<\theta_*$, at $c=c_*$.    
}\\

To justify assumption (i), consider first a general mixture of two non-interpenetrating 
phases 1 and 2 of varying mass fractions $c, 1-c$
and possibly varying specific volumes $\tau_1$ and  $\tau_2$, for which 
$$
\tau=c\tau_1+(1-c)\tau_2
$$ 
and 
$$\hat U=cU_1(\tau_1)+(1-c)U_2(\tau_2).$$ 
Then restrict attention to the case that phase 1 is perfectly incompressible and thus
does not store mechanical energy, 
$$
\tau_1=\hbox{const}\quad\hbox{ and }\quad U_1(\tau_1)=0.
$$
Supposing further 
that phase 2 is lighter than phase 1,
$$
\tau_2>\tau_1
$$
and that, as a simple prototypical example, its
energy has the form
$$
U_2(\tau_2)=-\log(\tau_2)
$$
leads to the stated form of $\hat U$ as a function of $\tau$ and $c$.
\par\medskip
Noticing that $\hat U_{\tau\tau}\hat U_{cc}-\hat U_{c\tau}^2=0$ and thus,
in the terminology of \cite{F3}, 
\be\label{sgnDelta}
\hbox{sgn}(\Delta(\tau,c))
=
\hbox{sgn}(W_{cc}(c,\theta)),
\ee
one readily sees that Theorem 1 follows exactly as in \cite{F3}.\footnote{It is 
actually easier, here and in other contexts, to work directly with the Gibbs potential 
$G$ associated with $U$. The role of $\Delta
=U_{\tau\tau} U_{cc}-U_{c\tau}^2$
is then played by the simpler quantity $G_{cc}$. In the present context, $\hat G_{cc}=0$ 
and Eq.\ \eqref{sgnDelta} reads 
$
\hbox{sgn}(G_{cc}(p,c))
=
\hbox{sgn}(W_{cc}(c,\theta))
$.}
We illustrate this by pointing out that 
by Lemma 1 of \cite{F3}, the proof amounts to
studying the level sets of 
$$
\Gamma^{\theta,\pi}(c,y)\equiv \hat G(P^\pi(c,y),c)+W(c,\theta)+\frac12 y^2,
$$
where $y$ corresponds to $c'$,
$$
\hat G(p,c)=(1-c)(1+\log p)+cp\tau_1
$$
is the Gibbs potential associated with $\hat U$ and $P^\pi(c,y)$ the unique positive root $p$ of 
$$
0=(p-\pi)(c\tau_1p+(1-c))
+{y^2p}.
$$
(The latter equation is Eq.\ (9) in  \cite{F3}.) 
The critical pressure is $p=p_*$,   
the unique solution $<1$ of
$$
G_c(p,c)=\tau_1p-\log p-1=0.
$$
For $(\theta,\pi)$ near $(\theta_*,p_*)$, the level landscape of $\Gamma^{\theta,\pi}$
undergoes a transition from one saddle (for $\theta>\theta_*$) to a saddle-maximum-saddle 
configuration (for $\theta<\theta_*$ and certain $\pi$). In the latter case, the two
saddles are at the same level and thus connected by two heteroclinic orbits (that together 
surround the maximum point) if $\pi$ assumes a unique value $\pi_*(\theta)$. As in \cite{F3}, 
Theorem 2 then follows from the transversality of the saddle-saddle
connections with respect to the parameter $\pi$. \par\medskip\bigskip\vskip -13cm
\includegraphics[height=200mm,width=185mm]{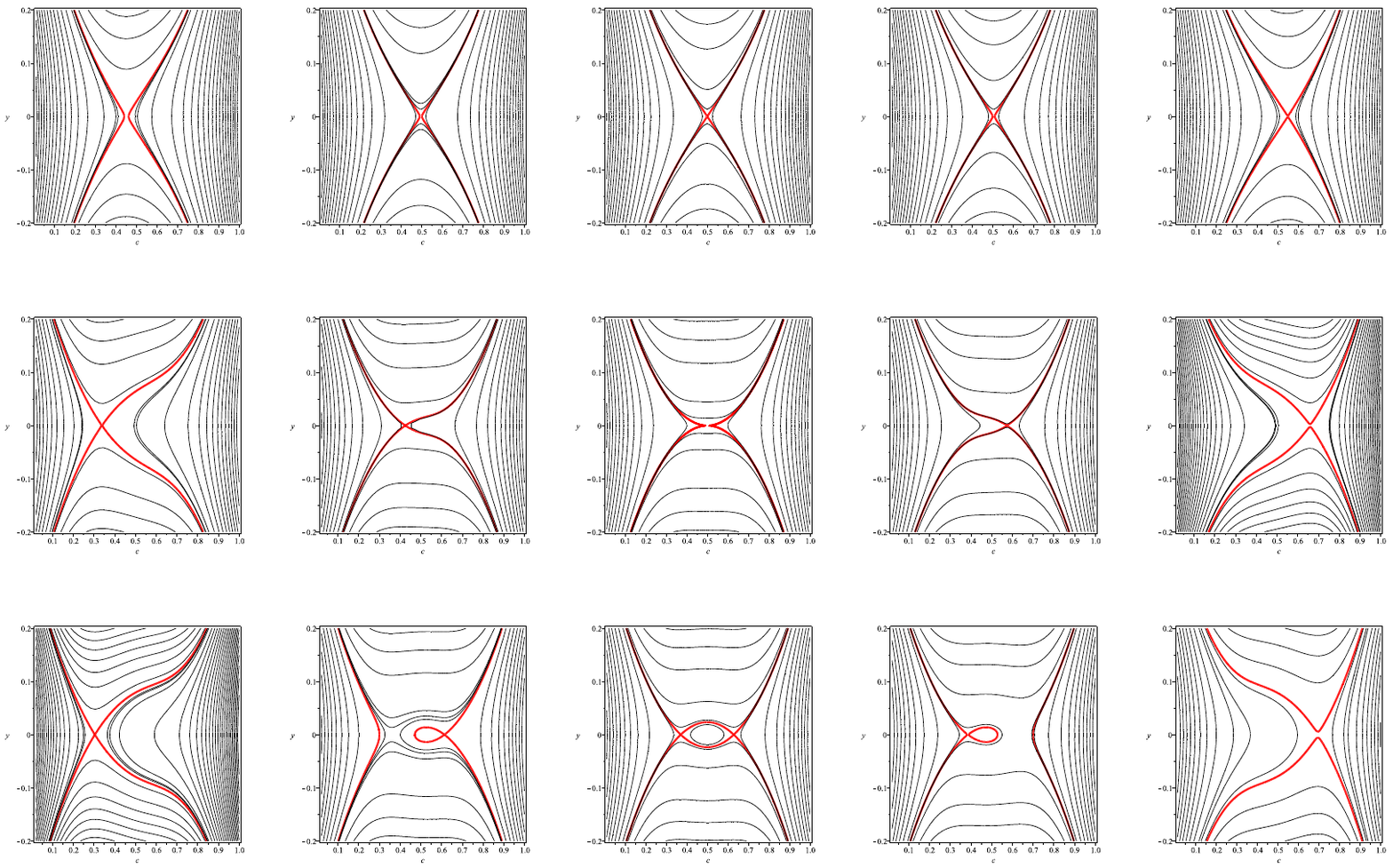}

\par
{\bf Figure} (by J. H\"owing){\bf :} 
Level lines of
$
\Gamma^{\theta,\pi}
$
for $\tau_1=0.5$ 
and $W(c,\theta)=(c-0.5)^4+(\theta-\theta_*)(c-0.5)^2$.\\
Top to bottom: $\theta-\theta_*= 0.16, 0.00, -0.08$. Left to right: $\pi-p_*=-0.010, -0.001, 0.000, 0.001, 0.010$.
\par\bigskip{\bf Remark.} {\it The choice of $-\log$ is exemplary. $U_2$ can be any 
function $f:(\tau_1,\infty)\to\R$ with $f'<0<f''$.}

\end{document}